\documentclass[12pt]{amsart}
\usepackage{amssymb,amsmath,amsthm,latexsym}
\usepackage{graphicx}
\newtheorem{theorem}{Theorem}

\newcommand{\be}{\begin{enumerate}}
\newcommand{\ee}{\end{enumerate}}

\newcommand{\beq}{\begin{equation}}
\newcommand{\eeq}{\end{equation}}
\newcommand{\beas}{\begin{eqnarray*}}
\newcommand{\eeas}{\end{eqnarray*}}
\newcommand{\bea}{\begin{eqnarray}}
\newcommand{\eea}{\end{eqnarray}}

\newcommand{\la}{\lambda}
\newcommand{\tw}{\tilde P}

\newcommand{\zz}{\mathbb{Z}}

\DeclareMathOperator{\diag}{diag}
\DeclareMathOperator{\SL}{SL}
\DeclareMathOperator{\GL}{GL}

\begin{document}

\title[Smith Normal Form of a Multivariate Matrix]{Smith Normal Form of a Multivariate Matrix Associated with
  Partitions
}

\author{Christine Bessenrodt}
\address{Institut f\"ur Algebra, Zahlentheorie und Diskrete Mathematik,
Leibniz Universit\"at Hannover,
D-30167 Hannover, Germany.}
\email{bessen@math.uni-hannover.de}

\author[Richard P. Stanley]{Richard P. Stanley${}^{\ast }$}
\address{Department of Mathematics, Massachusetts Institute
    of Technology, Cambridge, MA 02139, USA.}
\email{rstan@math.mit.edu}

\thanks{
\\
${}^*$This author's contribution is based upon work supported by the National
Science Foundation under Grant No.~DMS-1068625.}

\date{April 16, 2014 
}

\begin{abstract}
Consideration of a question of E.\ R.\ Berlekamp led Carlitz,
Roselle, and Scoville to give a combinatorial interpretation
of the entries of certain matrices of determinant~1  in terms of
lattice paths. Here we generalize this result by refining
the matrix entries to be multivariate polynomials, and by determining
not only the determinant but also the Smith normal form of these matrices.
A priori the Smith form need not exist
but its existence follows from the explicit computation.
It will be more convenient for us to state our results in
terms of partitions rather than lattice paths.
\end{abstract}

\maketitle

\medskip\bigskip

E. R. Berlekamp \cite{berl1,berl2} raised a question concerning
the entries of certain matrices of determinant 1.\ \ (Originally
Berlekamp was interested only in the entries modulo 2.) Carlitz,
Roselle, and Scoville \cite{c-r-s} gave a combinatorial interpretation
of the entries (over the integers, not just modulo 2) in terms of
lattice paths. Here we will generalize the result of Carlitz, Roselle,
and Scoville in two ways: (a) we will refine the matrix entries so
that they are multivariate polynomials, and (b) we compute
not just the determinant of these matrices, but more strongly their
Smith normal form (SNF). A priori our matrices need not have a Smith
normal form since they are not defined over a principal ideal domain,
but the existence of SNF will follow from its explicit computation. A
special case is a determinant of $q$-Catalan numbers.  It will be more
convenient for us to state our results in terms of partitions rather
than lattice paths.

Let $\lambda$ be a partition, identified with its Young
diagram regarded as a set of squares;
we  fix $\lambda$ for all that follows.
Adjoin to $\lambda$ a border
strip extending from the end of
the first row to the end of the first column of $\lambda$, yielding an
\emph{extended partition} $\lambda^*$. Let $(r,s)$ denote the square
in the $r$th row and $s$th column of $\lambda^*$. If
$(r,s)\in\lambda^*$, then
let $\lambda(r,s)$ be the partition whose diagram consists of all
squares $(u,v)$ of $\lambda$ satisfying $u\geq r$ and $v\geq s$. Thus
$\lambda(1,1)=\lambda$, while $\lambda(r,s)=\emptyset$ (the empty
partition) if $(r,s)\in\lambda^*\setminus\lambda$. Associate with the square
$(i,j)$ of $\lambda$ an
indeterminate~$x_{ij}$. Now for each square $(r,s)$ of~$\lambda^*$,
associate a polynomial $P_{rs}$ in the variables~$x_{ij}$, defined as
follows:
\beq
P_{rs} = \sum_{\mu\subseteq  \lambda(r,s)}\prod_{(i,j)\in\lambda(r,s)\setminus \mu} x_{ij},  \
  \label{eq:prsdef}
\eeq
where $\mu$ runs over all partitions contained in $\la (r,s)$.
In particular, if $(r,s)\in\lambda^*\setminus \lambda$ then $P_{rs}=1$.
Thus for $(r,s)\in\lambda$, $P_{rs}$ may be regarded as a generating
function for the squares
of all skew diagrams $\lambda(r,s)\setminus \mu$.  For instance, if
$\lambda=(3,2)$ and we set $x_{11}=a$, $x_{12}=b$, $x_{13}=c$,
$x_{21}=d$, and $x_{22}=e$, then Figure~\ref{fig1} shows the extended
diagram $\lambda^*$ with the polynomial $P_{rs}$ placed in the
square~$(r,s)$.

\begin{figure}
\setlength{\unitlength}{0.008in}
\begin{picture}(400,310)
\put(0,300){\line(1,0){400}}
\put(10,275){\scriptsize{$abcde+bcde$}}
\put(9,258){\scriptsize{$+bce+cde+$}}
\put(12,241){\scriptsize{$ce+de+c$}}
\put(20,224){\scriptsize{$+e+1$}}
\put(113,258){\scriptsize{$bce+ce+c$}}
\put(133,241){\scriptsize{$+e+1$}}
\put(234,248){\scriptsize{$c+1$}}
\put(345,248){\scriptsize{1}}
\put(0,210){\line(1,0){400}}
\put(18,158){\scriptsize{$de+e+1$}}
\put(133,158){\scriptsize{$e+1$}}
\put(245,158){\scriptsize{1}}
\put(345,158){\scriptsize{1}}
\put(0,120){\line(1,0){400}}
\put(0,30){\line(1,0){300}}
\put(45,68){\scriptsize{1}}
\put(145,68){\scriptsize{1}}
\put(245,68){\scriptsize{1}}
\multiput(0,300)(100,0){4}{\line(0,-1){270}}
\put(400,300){\line(0,-1){180}}
\end{picture}
\caption{The polynomials $P_{rs}$ for $\lambda=(3,2)$}
\label{fig1}
\end{figure}

\medskip

Write
  $$ A_{rs}=\prod_{(i,j)\in\lambda(r,s)} x_{ij}. $$
Note that $A_{rs}$ is simply the leading term of $P_{rs}$.
Thus for $\lambda=(3,2)$ as in Figure~\ref{fig1} we have
$A_{11}=abcde, A_{12}=bce$, $A_{13}=c$, $A_{21}=de$, and
$A_{22}=e$.

For each square $(i,j)\in\lambda^*$ there will be a unique subset of
the squares of $\lambda^*$ forming an $m\times m$ square $S(i,j)$ for
some $m\geq 1$, such that the upper left-hand corner of $S(i,j)$ is
$(i,j)$, and the lower right-hand corner of $S(i,j)$ lies in
$\lambda^*\setminus \lambda$. In fact, if $\rho_{ij}$ denotes the \emph{rank}
of $\lambda(i,j)$ (the number of squares on the main diagonal, or
equivalently, the largest $k$ for which $\lambda(i,j)_k\geq k$), then
$m=\rho_{ij}+1$. Let $M(i,j)$ denote the matrix obtained by inserting
in each square $(r,s)$ of $S(i,j)$ the polynomial $P_{rs}$. For
instance, for the partition $\lambda=(3,2)$ of Figure~\ref{fig1},
the matrix $M(1,1)$ is given by
 $$ M(1,1) = \left[ \begin{array}{ccc}
  P_{11} & bce+ce+c+e+1 & c+1\\ de+e+1 & e+1 & 1\\
  1 & 1 & 1 \end{array} \right], $$
where $P_{11}=abcde+bcde+bce+cde+ce+de+c+e+1$. Note that for this
example we have
  $$ \det M(1,1)= A_{11}A_{22}A_{33}=abcde\cdot e\cdot 1=abcde^2. $$

If $R$ is a commutative ring (with identity), and $M$ an $m\times n$
matrix over~$R$, then we say that $M$ has a \emph{Smith normal form}
(SNF) over~$R$ if there exist matrices $P\in\GL(m,R)$
(the set of $m\times m$ matrices over~$R$ which have an inverse whose
entries also lie in~$R$, so that $\det P$ is a unit in~$R$), $Q\in\GL(n,R)$,
such that $PMQ$ has the form (w.l.o.g., here $m\leq n$, the other case is dual)
$$ \begin{array}{rcl}
PMQ & = &
\left[ \begin{array}{ccccc}
 {\bf 0} & d_1 d_2 \cdots d_m &  &  & \\
 {\bf 0} & & d_1 d_2 \cdots d_{m-1} & {\bf 0} &   \\
  \vdots &  {\bf 0} &\ddots &  \\
 {\bf 0} & & &  & d_1
\end{array}
\right]
\\[10pt]
&=& ({\bf 0}, \mathrm{diag}(d_1 d_2 \cdots d_m, d_1 d_2 \cdots d_{m-1},\dots,
d_1))\:,
\end{array} $$
where each $d_i\in R$. If $R$ is an integral domain and $M$
has an SNF, then it is unique up to multiplication of the diagonal
entries by units. If $R$ is a principal ideal domain then the SNF
always exists, but not over more general rings. We will be working
over the polynomial ring
 \beq R=\zz[x_{ij}\colon (i,j)\in\lambda]. \label{eq:rdef} \eeq
Our main result asserts that $M(i,j)$ has a Smith normal form over~$R$,
which we describe explicitly. In particular, the entries on the
main diagonal are monomials.
It is stated below for $M(1,1)$, but it applies to any
$M(i,j)$ by replacing $\lambda$ with $\lambda(i,j)$.
Note also that the transforming matrices
are particularly nice as they are triangular matrices
with 1's on the diagonal.

\begin{theorem} \label{thm:snf}
There are an upper unitriangular matrix $P$
and a lower unitriangular matrix $Q$ in $\SL(\rho+1,R)$
such that
   $$ P\cdot M(1,1)\cdot Q =
      \mathrm{diag}(A_{11},A_{22},\dots,A_{\rho+1,\rho+1}). $$
In particular, $\det M(1,1)=A_{11}A_{22}\cdots A_{\rho\rho}$ (since
$A_{\rho+1,\rho+1}=1$).
\end{theorem}

For instance, in the example of Figure~\ref{fig1} we have
$$P \cdot M(1,1)\cdot Q=\mathrm{diag}(abcde,e,1)\:.$$

\medskip

We will give two proofs for Theorem~\ref{thm:snf}.

The main tool used for the first proof
is a recurrence for the polynomials~$P_{rs}$. To state
this recurrence we need some definitions. Let
$\rho=\mathrm{rank}(\lambda)$, the size of the main diagonal of
$\lambda$. For $1\leq i\leq \rho$, define
the rectangular array
  $$ R_i = \left[ \begin{array}{cccc}
     (1,2) & (1,3) & \cdots & (1,\lambda_i-i+1)\\
     (2,3) & (2,4) & \cdots & (2,\lambda_i-i+2)\\
           &       & \vdots\\
      (i,i+1) & (i,i+2) & \cdots & (i,\lambda_i) \end{array}
  \right] \:.$$
Note that if $\lambda_\rho=\rho$, then $R_\rho$ has no columns, i.e.,
$R_\rho=\emptyset$.
Let $X_i$ denote the set of all subarrays of $R_i$ whose shapes form a
vertical reflection of a Young diagram, justified into the upper
right-hand corner of $R_i$. Define
   $$ \Omega_i =\sum_{\alpha\in X_i}\prod_{(a,b)\in\alpha} x_{ab}. $$
For instance, if $\lambda_2=4$, then
   $$ R_2 = \left[ \begin{array}{cc} (1,2) & (1,3)\\ (2,3) & (2,4)
     \end{array} \right], $$
so
   $$ \Omega_2 = 1+x_{13}+x_{12}x_{13}+x_{13}x_{24}
     +x_{12}x_{13}x_{24}+x_{12}x_{13}x_{23}x_{24}. $$
In general, $R_i$ will have $i$ rows and $\lambda_i-i$ columns, so
$\Omega_i$ will have $\binom{\lambda_i}{i}$ terms. We also set
$\Omega_0=1$, which is consistent with regarding $R_0$ to be an empty
array.

Next set $S_0=\emptyset$ and define for $1\leq i\leq \rho$,
  $$ S_i = \{ (a,b)\in\lambda\colon 1\leq a\leq i,\ \
    \lambda_i-i+a< b\leq \lambda_a\}. $$
In particular $S_1=\emptyset$. When $\lambda_\rho=\rho$,
then $S_\rho$ consists of all squares strictly to the right of the
main diagonal of $\lambda$;
otherwise, $S_i$ consists of those squares of $\lambda$ that are
in the same row and to the right of all squares appearing as an entry
of $R_i$.
\\
Set
   $$ \tau_i = \Omega_i\cdot \prod_{(a,b)\in S_i} x_{ab}, $$
where as usual an empty product is equal to 1. In particular,
$\tau_0=1$.

We can now state the recurrence relation for $P_{rs}$.

\begin{theorem} \label{thm:prs}
Let $2\leq j\leq \rho+1=\mathrm{rank}(\lambda)+1$. Then
  $$ \tau_0  P_{1j}-\tau_1 P_{2j} + \cdots +
     (-1)^\rho\tau_\rho P_{\rho+1, j}=0.
  $$
Moreover, for $j=1$ we have
    $$ \tau_0    P_{11}-\tau_1 P_{21} + \cdots +
     (-1)^\rho\tau_\rho P_{\rho+1, 1}=A_{11}.
 $$
\end{theorem}

Before presenting the proof, we first give a couple of examples. Let
$\lambda=(3,2)$, with  $x_{11}=a$, $x_{12}=b$,
$x_{13}=c$, $x_{21}=d$, and $x_{22}=e$ as in Figure~\ref{fig1}. For
$j=1$ we obtain the identity
  $$ (1+c+e+ce+de+bce+cde+bcde+abcde)-(1+c+bc)(1+e+de)+bc $$
\vspace{-2em} \hspace{-3em}
  $$  = abcde, $$
For $j=2$ we have
   $$ (1+c+e+ce+bce)-(1+c+bc)(1+e)+bc = 0, $$
while for $j=3$ we get
   $$ (1+c)-(1+c+bc)+bc=0. $$
For a further example, let $\lambda=(5,4,1)$, with the variables
$x_{ij}$ replaced by the letters $a,b,\dots,j$ as shown in
Figure~\ref{fig2}. \\
For $j=1$ we get
  $$ P_{11} -(1+e+de+cde+bcde)(1+i+j+ij+hi+hij+ghi+ghij+fghij)  $$
  $$ +de(1+c+bc+ci+bci+bchi)(1+j) = abcdefghij, $$
where $P_{11}=1+e+i+j+\cdots+abcdefghij$, a polynomial with 34
terms. For $j=2$ we get
 $$ (1+e+i+ei+hi+dei+ehi+ghi+dehi+eghi+cdehi+deghi+cdeghi+ bcdeghi) $$
 $$ - (1+e+de+cde+bcde)(1+i+hi+ghi)+de(1+c+bc+ci+bci+bchi)=0. $$
For $j=3$ we have
   $$ (1+e+i+ei+hi+dei+ehi+dehi+cdehi) -
      (1+e+de+cde+bcde)(1+i+hi)  $$
   $$ +de(1+c+bc+ci+bci+bchi)=0. $$

\begin{figure}
\setlength{\unitlength}{0.004in}%
\begin{picture}(500,300)
\multiput(0,300)(0,-100){2}{\line(1,0){500}}
\put(45,240){{$a$}}
\put(145,240){{$b$}}
\put(245,240){{$c$}}
\put(345,240){{$d$}}
\put(445,240){{$e$}}
\put(45,140){{$f$}}
\put(145,140){{$g$}}
\put(245,140){{$h$}}
\put(345,140){{$i$}}
\put(0,100){\line(1,0){400}}
\put(45,40){{$j$}}
\multiput(0,300)(100,0){2}{\line(0,-1){300}}
\multiput(200,300)(100,0){3}{\line(0,-1){200}}
\put(500,300){\line(0,-1){100}}
\put(0,0){\line(1,0){100}}
\end{picture}
\caption{The variables for $\lambda=(5,4,1)$}
\label{fig2}
\end{figure}

\textbf{Proof of Theorem~\ref{thm:prs}.} First suppose that $j\geq
2$. We will prove the result in the form
  $$ \tau_0  P_{1j}= \tau_1 P_{2j} - \cdots +
     (-1)^{\rho-1}\tau_\rho P_{\rho+1,j}
  $$
by an Inclusion-Exclusion argument. Since $\tau_0=1$, the
left-hand side is the generating function for all skew diagrams
$\lambda(1,j)\setminus \mu$, as defined by equation~\eqref{eq:prsdef}. If we
take a skew diagram $\lambda(2,j)\setminus \sigma$ and append to it some
element of~$X_1$ (that is, some squares on the first row forming a
contiguous strip up to the last square $(1,\lambda_1)$), then we will
include every skew diagram $\lambda(1,j)\setminus \mu$. However, some
additional diagrams $\delta$ will also be included. These will have
the property that the first row begins strictly to the left of the
second. We obtain the first two rows of such a diagram $\delta$ by
choosing an element of~$X_2$ and adjoining to it the set~$S_2$. The
remainder of the diagram $\delta$ is a skew
shape~$\lambda(3,j)\setminus \zeta$. Thus we cancel out the unwanted terms of
$\tau_1 P_{2j}$ by subtracting $\tau_2 P_{3j}$.
However, the product $\tau_2 P_{3j}$ has some
additional terms that need to be added back in. These terms will
correspond to diagrams $\eta$ with the property that the first row
begins strictly to the left of the second, and the second begins
strictly to the left of the third.  We obtain the first three rows of
such a diagram $\eta$ by choosing an element of~$X_3$ and adjoining to
it the set $S_3$. The remainder of the diagram $\eta$ is a
skew shape~$\lambda(4,j)\setminus \xi$. Thus we cancel out the unwanted terms of
$\tau_2 P_{3j}$ by adding $\tau_3 P_{4j}$. This
Inclusion-Exclusion process will come to end when we reach the term
$\tau_\rho P_{\rho+1,j}$, since we cannot have $\rho+1$
rows, each strictly to the left of the one below. This proves the
theorem for $j\geq 2$.

When $j=1$, the Inclusion-Exclusion process works exactly as before,
except that the term $A_{11}$ is never cancelled from
$\tau_0 P_{11} = P_{11}$. Hence the theorem is also true for $j=1$.
$\ \ \Box$

\medskip

With this result at hand, we can now embark on the \textbf{proof of
  Theorem \ref{thm:snf}.}
This is done by induction on $\rho$, the
result being trivial for $\rho=0$ (so $\lambda=\emptyset$).
Assume the assertion holds for
partitions of rank less than $\rho$, and let rank$(\lambda)=\rho$. For
each $1\leq i\leq \rho$,
multiply row $i+1$ of $M(1,1)$ by $(-1)^i\tau_i$ and add it to
the first row. By Theorem~\ref{thm:prs} we get a matrix $M'$ whose
first row is
$[A_{11},0,0,\dots,0]$. Now by symmetry we can perform the analogous
operations on the \emph{columns} of $M'$. We then get a
matrix in the block
diagonal form
$\left[ \begin{array}{cc} A_{11} & 0\\ 0 & M(2,2) \end{array} \right]$.
The row and column operations that we have performed are
equivalent to computing $P'MQ'$ for upper and lower unitriangular
matrices $P',Q'\in \SL(\rho+1,R)$, respectively.
The proof now follows by induction. $\ \ \Box$

\textsc{Note.}  The determinant above can also
  easily be evaluated by the Lindstr\"om-Wilf-Gessel-Viennot method of
  nonintersecting lattice paths, but it seems impossible to extend
  this method to a computation of SNF.
\medskip

We now come to the second approach towards the SNF
which does not use Theorem~\ref{thm:prs}.
Indeed, we will prove the more general version below where
the weight matrix is not necessarily square;
while this is not expounded here, the previous proof may also
be extended easily to any rectangular
subarray (regarded as a matrix) of $\lambda^\ast$ whose lower-right
hand corner lies in $\lambda^\ast\setminus \lambda$.
The inductive proof below will not involve Inclusion-Exclusion arguments;
again, suitable transformation matrices are computed
explicitly stepwise along the way.

Given our partition $\la$,
let $F$ be a rectangle of size~$d\times e$ in $\la^*$,
with top left corner at $(1,1)$, such
that its corner $(d,e)$ is a
square in the added border strip $\la^*\setminus\la$;
thus $A_{de}=1$ and $P_{de}=1$.
We denote the corresponding matrix of weights by
$$W_F=(P_{ij})_{(i,j)\in F}\:.$$

\begin{theorem}\label{thm:snf-rect}
Let $F$ be a rectangle in $\la^*$ as above, of size~$d\times e$;
assume $d\leq e$ (the other case is dual).
Then there are an upper unitriangular matrix $P \in \SL(d,R)$
and a lower unitriangular matrix $Q \in \SL(e,R)$
such that
$$ P\cdot W_F \cdot Q =
({\bf 0}, \diag(A_{1,1+e-d},A_{2,2+e-d}\ldots, A_{d,e}))
\:. $$
In particular, when $\la$ is of rank $\rho$ and
$F$ is the Durfee square in $\la^*$, we have the result
in Theorem~\ref{thm:snf} for $W_F=M(1,1)$.
\end{theorem}

\textbf{Proof.}
We note that the claim clearly holds
for $1\times e$ rectangles as then $A_{1,e}=1=P_{1,e}$.
We use induction on the size of $\la$.
For  $\la=\emptyset$, $\la^*=(1)$ and $F$ can only be a
$1\times 1$ rectangle.
We now assume that the result holds for all rectangles as above
in partitions of~$n$.
We take a partition~$\la'$ of $n+1$ and a rectangle
$F'$ with its corner on the rim $\la'^*\setminus \la'$,
where we may assume that $F'$ has at least two rows;
we will produce the required transformation matrices
inductively along the way.

First we assume that we can remove a square $s=(a,b)$
from $\la'$ and obtain a partition $\la=\la'\setminus{s}$
such that $F'\subseteq \la^* $.
By induction, we thus have the result for $\la$ and
$F=F'\subset \la^*$; say this is
a rectangle with corner $(d,e)$, $d\le e$.
Then $a<d$ and $b\ge e$, or $a\ge d$ and $b<e$.
We discuss the first case, the other case is analogous.
Set $z=x_{ab}$.

Let $t\in \la\subset \la'$; denote the weights of $t$
with respect to~$\la$ by $A_t$, $P_t$, and
with respect to~$\la'$ by $A'_t$, $P'_t$.
Let $W=W_F=(P_{ij})_{(i,j)\in F}$
and $W'=(P'_{ij})_{(i,j)\in F}$
be the corresponding weight matrices.
We clearly have
$$A'_{j,j+e-d} =
\left\{
\begin{array}{rl}
z A_{j,j+e-d} & \text{ for } 1\le j\le a \\
  A_{j,j+e-d} & \text{ for } a < j \le d
\end{array}
\right. \:.
$$
When we compute the weight of  $t=(i,j)\in F$ with $i\le a$
with respect to~$\la'$,
we get two types of contributions to~$P_{ij}'$.
For a partition
$\mu \subseteq \la'(i,j)$
with $s\not\in \mu$, i.e., $\mu\subseteq \la(i,j)$,
the corresponding weight summand in $P_{ij}$ is multiplied by~$z$,
and hence the total contribution
from all such $\mu$ is exactly $zP_{ij}$.
On the other hand, if the partition
$\mu \subseteq \la'(i,j)$ contains $s$,
then it contains the full rectangle $\mathcal R$
spanned by the corners $t$ and $s$, that is,
running over rows $i$ to $a$ and columns $j$ to $b$;
$\mu$ arises from $\mathcal R$ by gluing suitable
(possibly empty) partitions $\alpha$, $\beta$
to its right and bottom side, respectively, and
$$\la'(i,j) \setminus \mu = (\la(i,b+1)\setminus \alpha) \cup
(\la(a+1,j)\setminus \beta ) \:.$$
Summing the terms for all such $\mu$, we get the contribution
$P_{i,b+1}\cdot P_{a+1,j}$.
Clearly, when  $i>a$, then the square $s$
has no effect on the weight of~$t$.
Hence
$$
P_{ij}'
 =
\left\{
\begin{array}{cl}
z P_{ij} +  P_{i,b+1}\cdot P_{a+1,j} & \text{ for } 1\le i\le a \\
  P_{ij} & \text{ for } a < i\le d
\end{array}
\right. \:.
$$
We now transform $W'$.
Multiplying the $(a+1)$-th row of $W'$ by $P_{i,b+1}$ and subtracting
this from row $i$, for all $i\le a$, corresponds to multiplication
from the left with an upper unitriangular matrix in $\SL(d,R)$
and gives
$$W'_1=
\left[ \begin{array}{ccc}
zP_{11} & \cdots &  zP_{1 e}\\
\vdots & \cdots  & \vdots \\
zP_{a1} &  \cdots & zP_{a e }\\
P_{a+1,1} &  \cdots & P_{a+1,e }\\
\vdots &  \cdots & \vdots \\
P_{d 1}  & \cdots & P_{d e }
\end{array} \right]
\:.
$$
By induction, we know that there are upper and lower unitriangular matrices
$U=(u_{ij})_{1\le i,j\le d },V=(v_{ij})_{1\le i,j\le e }$, respectively,
defined over~$R$, such that
$UWV=( {\bf 0},\diag(A_{1,1+e-d},\ldots,A_{d,e }))$.
We then define an upper unitriangular
matrix $U'=(u_{ij}')_{1\le i,j\le d }\in \SL(d,R)$
by setting
$$u'_{ij} =
\left\{
\begin{array}{rl}
z u_{ij} & \text{ for } i\le a <j \\
  u_{ij} & \text{ otherwise }
\end{array} \:.
\right.
$$
With $UW = (\tilde P_{ij})$, we then have
$$U'W_1'= 
\left[ \begin{array}{ccc}
z\tw_{11} & \cdots & z\tw_{1 e }\\
\vdots & \cdots & \vdots \\
z\tw_{a1} & \cdots & z\tw_{a e }\\
\tw_{a+1,1} & \cdots & \tw_{a+1,e }\\
\vdots & \cdots & \vdots \\
\tw_{d 1} & \cdots & \tw_{d e }
\end{array} \right],
$$
and hence we obtain the desired form via
$$U'W'_1V=( {\bf 0},\diag(zA_{1,1+e-d},\ldots,zA_{a,a+e-d},A_{a+1,a+1+e-d},\ldots, A_{d e }))\:.$$

Next we deal with the case where we cannot remove a square
from $\la'$ such that the rectangle $F'$ is still contained
in the extension of the smaller partition~$\la$ of~$n$.
This is exactly the case when $\la'$ is a rectangle,
with corner square $s=(d,e)$ (say),
and $F'=\la'^*$ is the rectangle with its corner at $(d+1,e+1)$.
Then $s$ is the only square that can be removed from~$\la'$;
for $\la=\la'\setminus s$ we have $\la^*=\la'^*\setminus (d+1,e+1)$.
We now use the induction hypothesis for the partition $\la$ of $n$
and the rectangle $F=\la'\subset \la'^*=F'$.

We keep the notation for the weights
of a square $t=(i,j)$ with respect to~$\la$, $\la'$,
and set $z=x_s=x_{de}$.
Clearly, we have for the monomial weights to~$\la'$:
$$A'_{j,j+e-d} =
\left\{
\begin{array}{cl}
z A_{j,j+e-d} & \text{ for } j\le d  \\
  1 & \text{ for } j = d +1
\end{array}
\right. \:.
$$
Now we consider how to compute
the matrix $W'=(P'_{ij})_{(i,j)\in F'}$
from the values $P_{ij}$, $(i,j)\in F$.
Arguing analogously as before,
we obtain
$$
P_{ij}'
 =
\left\{
\begin{array}{cl}
z P_{ij} +  1 & \text{ for } 1\le i\le d, 1\le j \le e  \\
   1& \text{ for } i\le d+1  \text { and } j=e +1 \\
   1 & \text{ for } i=d +1 \text{ and } j\le e+1
\end{array}
\right. \:.
$$
As a first simplification on $W'$,
we subtract the $(d +1)$-th row of $W'$  from row $i$, for all $i\le d $ (corresponding to a multiplication
from the left with a suitable upper unitriangular matrix), and we obtain
the matrix
$$W'_1= 
\left[ \begin{array}{cccc}
zP_{11} & \ldots & zP_{1 e } & 0\\
\vdots & \ldots & \vdots & \vdots \\
zP_{d 1} & \ldots & zP_{d e } & 0 \\
1 & \ldots & 1 & 1
\end{array} \right]\:.
$$
Subtracting the last column from column $j$, for all $j\le e $,
transforms this (via
postmultiplication with a lower unitriangular matrix as required)
into
$W'_2= 
\left[ \begin{array}{cc}
zW & {\bf 0}\\
{\bf 0} & 1
\end{array} \right]$.
By induction we have an upper and a lower unitriangular matrix
$U\in \SL(d,R)$, $V\in \SL(e,R)$, respectively, with
$UWV=( {\bf 0},\diag(A_{1,1+e-d}, \ldots, A_{d e }))$.
Then
$$
\left[ \begin{array}{cc}
U & {\bf 0}\\
{\bf 0} & 1
\end{array} \right]
\left[ \begin{array}{cc}
zW & {\bf 0}\\
{\bf 0} & 1
\end{array} \right]
\left[ \begin{array}{cc}
V & {\bf 0}\\
{\bf 0} & 1
\end{array} \right]
=
\left[ \begin{array}{ccc}
 {\bf 0} & \diag(zA_{1,1+e-d}, \ldots, zA_{d e}) & {\bf 0} \\
{\bf 0} & {\bf 0}& 1
\end{array} \right]\:,  $$
and we have the assertion as claimed. $\Box$

\medskip

We conclude with an interesting special case, namely, when
$\lambda$ is the ``staircase''
$(n-1,n-2,\dots,1)$ and each $x_{ij}=q$, we get that $P_{ij}$ is the
$q$-Catalan number that is denoted $\tilde{C}_{n+2-i-j}(q)$ by
F\"urlinger and Hofbauer \cite{f-h},\cite[Exer.~6.34(a)]{ec2}. For
instance, $\tilde{C}_3(q)=1+2q+q^2+q^3$. Theorem~\ref{thm:snf} gives
that the matrix
  $$ M_{2m-1}=[\tilde{C}_{2m+1-i-j}(q)]_{i,j=1}^m $$
has SNF
$\mathrm{diag}(q^{\binom{2m-1}{2}},q^{\binom{2m-3}{2}},q^{\binom{2m-5}{2}},
\dots,q^3,1)$, and the matrix
  $$ M_{2m} =[\tilde{C}_{2m+2-i-j}(q)]_{i,j=1}^{m+1} $$
has SNF $\mathrm{diag}(q^{\binom{2m}{2}},q^{\binom{2m-2}{2}},
q^{\binom{2m-4}{2}},\dots,q,1)$. The determinants of the matrices
$M_n$ were already known (e.g., \cite{cigler1},\cite[p.~7]{cigler2}),
but their Smith normal form is new.

\end{document}